\documentclass[graybox, envcountchap]{svmult}

\newcommand{\mycomment}[1]{%
\ifthenelse{\isodd{\value{page}}}{%
\normalmarginpar%
\marginpar{\tiny {#1}}%
}{%
\reversemarginpar%
\marginpar{\tiny {#1}}%
}}%

\usepackage{graphicx}           
\usepackage{makeidx}            
\usepackage{multicol}           
\usepackage[bottom]{footmisc}   

\usepackage[utf8]{inputenc}
\usepackage[T1]{fontenc}
\usepackage{float}
\usepackage{wrapfig}
\usepackage{amsmath}
\usepackage{amsopn}
\usepackage{amsfonts}
\usepackage{hyperref}
\usepackage{lmodern}
\usepackage[english]{babel}
\usepackage[font=scriptsize,labelfont=bf]{caption}
\usepackage{booktabs}  
\usepackage{enumitem}  
\usepackage{tabularx}
\usepackage[normalem]{ulem}

\usepackage{enumitem}

\usepackage[ruled,vlined]{algorithm2e}
\SetAlgorithmName{Algorithm}{Algorithm}{List of Algorithms}

\usepackage{epstopdf}
\ifpdf
  \DeclareGraphicsExtensions{.eps,.pdf,.png,.jpg}
\else
  \DeclareGraphicsExtensions{.eps}
\fi

\newcommand{\N}{\mathbb{N}}
\newcommand{\grad}{{\nabla}}
\newcommand{\utilde}{\tilde{u}}
\newcommand{\calL}{\mathcal{L}}
\newcommand{\lambdak}[1][k]{\lambda^{(#1)}}
\newcommand{\alphak}[1][k]{\alpha_{#1}}

\spnewtheorem{assumption}{Assumption}[section]{\bf}{\it}
\newcommand{\TheTitle}{Combining Nonlinear FETI-DP Methods and Quasi-Newton Methods using an SQP Approach}
\newcommand{\TheAuthors}{S.\ Köhler, and O.\ Rheinbach}

\ifpdf
\hypersetup{
  pdftitle={\TheTitle},
  pdfauthor={\TheAuthors}
}
\fi

\definecolor{darkorange}{HTML}{B45F04}
\definecolor{royalpurple}{rgb}{0.47, 0.32, 0.66}

\begin{document}

\title*{{\TheTitle}}
\titlerunning{Nonlinear FETI-DP Methods and Quasi-Newton Methods}
\author{S.\ Köhler\orcidID{0000-0003-1015-8736}, and O.\ Rheinbach\orcidID{0000-0002-9310-8533}}
\institute{Stephan Köhler \and Oliver Rheinbach {\at} Technische Universität Bergakademie Freiberg, Akademiestr.~6, 09596~Freiberg \email{ stephan.koehler@math.tu-freiberg.de, oliver.rheinbach@math.tu-freiberg.de}
}

\numberwithin{theorem}{section}
\numberwithin{equation}{section}

\maketitle

\abstract*{
  The combination of nonlinear FETI-DP (Dual Primal Finite Element Tearing and Interconnecting) and Quasi-Newton methods using a sequential quadratic programming (SQP) approach is considered.   Nonlinear FETI-DP methods are parallel iterative solution methods for nonlinear finite element problems, based on divide and conquer, using Lagrange multipliers.  In the method, we use Quasi-Newton approximations of Hessian for the quadratic programs, where the initial approximation uses the exact Hessian.  To accelerate the convergence, we recompute the initial Hessian and restart the Quasi-Newton approximation.  We provide numerical experiments using homogeneous model problems from nonlinear structural mechanics.
}

\section{Nonlinear FETI-DP}\label{sec:feti-dp-basics}

Nonlinear FETI-DP methods~\cite{klawonn2017nonlinear} are nonlinear generalizations of standard linear dual-primal Finite Element Tearing and Interconnecting domain decomposition methods~\cite{toselli:2005:ddm} based on a divide-and-conquer approach. We consider an unconstrained minimization problem for some objective $J$ and decompose it following the FETI-DP strategy.  This unconstrained problem is transformed into a constrained one by tearing the computational domain at the interface into subdomains, using subassembly for the primal degrees of freedom (dofs) of the subdomain boundaries, and equality constraints for the remaining dofs of the subdomain boundaries, the so called dual variables.  The resulting problem can be written as
\begin{equation}
\textstyle  \label{eq:minimization-problem}
  \min_{\utilde} \widetilde{J}(\utilde) \quad \text{subject to (s.t.)} \quad B\utilde = 0,
\end{equation}
where the constraint $B\utilde = 0$ enforces continuity of the dual variables across subdomain boundaries; here,
$\utilde := [u_{B}^{(1)},\ldots,u_{B}^{(N)},\utilde_{\Pi}]^{T}$. The subscript $B$ refers to the union of the inner variables of the subdomains and dual ones, $\widetilde{J}(\utilde) := \sum_{i=1}^{N}J^{(i)}(u_{B}^{(i)},R_{\Pi}^{(i)}\utilde_{\Pi})$, $J^{(i)}$ is the local objective of the $i$-th subdomain, and $R_{\Pi}^{(i)}$ is the assembly operator of the primal variables as in linear FETI-DP methods~\cite{toselli:2005:ddm}.  The Lagrange function for~\eqref{eq:minimization-problem} is $\calL(\utilde,\lambda) = \widetilde{J}(\utilde) + \lambda^{T}\,B\utilde$.  The (nonlinear) saddle point problem of the first-order necessary optimality condition is
\begin{equation}
  \label{eq:global-saddle-point-problem}
  \begin{array}{rcr rcl}
    \grad_{\utilde}\calL(\utilde,\lambda) &=& \grad \widetilde{J}(\utilde) + B^T \lambda &=& \tilde{f} ,\\
    \grad_{\lambda}\calL(\utilde,\lambda) &=& B \tilde{u} &=& 0.
  \end{array}
\end{equation}
The coarse problem of the method is nonlinear and is obtained from finite element subassembly $\nabla \widetilde{J}(\utilde):=R_{\Pi}^{T}\nabla J(R_{\Pi}\tilde{u})$, as in linear FETI-DP methods~\cite{toselli:2005:ddm}.

\section{A Sequential Quadratic Programming Approach}\label{sec:a-sequential-quadratic-programming-approach}

We consider the sequential quadratic programming (SQP) approach for the constrained optimization problem~\eqref{eq:minimization-problem}, see, e.g.,~\cite{nocedal2006numerical, ulbrich2012nichtlineare, bertsekas1982constrained}.  The original constrained problem is replaced by a sequence of quadratic problems of the form
\begin{equation}
  \label{eq:sqp-subproblem-1}
  \begin{aligned}
    &\min_{d}  \grad_{\utilde}{\calL(\utilde^{(k)},\lambdak)}^{T}d + \frac{1}{2}d^{T}\grad^{2}_{\utilde\utilde}\calL(\utilde^{(k)},\lambdak)d \\
    \text{s.t.}\qquad & B(\utilde^{(k)} + d) = 0,
  \end{aligned}
\end{equation}
where $(\utilde^{(k)}, \lambdak)$ is the current iterate.  The next iterate is obtained by
\begin{equation}
  \label{eq:sqp-update-iterate}
  \begin{aligned}
    \utilde^{(k+1)} = \utilde^{(k)} + \alphak\delta\utilde^{(k)},\qquad
    \lambda^{(k+1)} = \delta\lambdak
  \end{aligned}
\end{equation}
where $(\delta\utilde^{(k)},\delta\lambdak)$ is the Karush-Kuhn-Tucker (KKT) point of the SQP subproblem~\eqref{eq:sqp-subproblem-1} and $\alphak$ is a suitable step length.  To determine the step length $\alphak$, we use the exact nondifferentiable penalty function
\begin{equation}
  \label{eq:penalty-function-def}
  P_{1}(\utilde;\mu) = \widetilde{J}(\utilde) + \mu\|B\utilde\|_{1},
\end{equation}
where $\mu$ is the penalty parameter.  It can be shown that for $P_{1}$ the one-sided directional derivative $DP_{1}(\utilde;d,\mu)$ exists for all points $\utilde$ and in all directions $d$.  Hence, we can use backtracking and the Armijo condition~\cite{armijo1966minimization}
\begin{equation}
  \label{eq:penalty-function-armijo}
  P_{1}(\utilde^{(k)}+\alphak\delta\utilde^{(k)};\mu) - P_{1}(\utilde^{(k)};\mu) \leq \alphak DP_{1}(\utilde;d,\mu)
\end{equation}
to compute the step length $\alphak$, see, e.g.,~\cite{nocedal2006numerical}.
Notice, for the directional derivative $DP_{1}$ it holds that
\begin{equation}
  \label{eq:penalty-function-dir-der}
  \begin{aligned}
    DP_{1}(\utilde;\delta\utilde^{(k)},\mu) =& \grad\widetilde{J}(\utilde^{(k)}) - \mu\|B\utilde^{(k)}\|_{1}\\
    \leq& -\left.\delta\utilde^{(k)}\right.^{T}\grad^{2}_{\utilde\utilde}\calL(\utilde^{(k)},\lambdak)\delta\utilde^{(k)} 
    -(\mu-\|\delta\lambdak\|_{\infty})B\utilde^{(k)},
  \end{aligned}
\end{equation}
where $(\delta\utilde^{(k)},\delta\lambdak)$ is the KKT point of~\eqref{eq:sqp-subproblem-1}, see, e.g.,~\cite{nocedal2006numerical}.  Therefore, if $\grad^{2}_{\utilde\utilde}\calL(\utilde^{(k)},\lambdak)$ is positive definite on $\ker(B)$, we obtain a descend direction for $P_{1}$ by the solution of~\eqref{eq:sqp-subproblem-1}.  A globalized SQP algorithm is presented in Fig.~\ref{fig:globalized-sqp}.  It holds, that under sufficient conditions all accumulation points of the sequence of iterates ${\left(\utilde^{(k)}\right)}_{k\in\N}$, generated by the algorithm in Fig.~\ref{fig:globalized-sqp}, are KKT points of the original problem~\eqref{eq:minimization-problem}, see, e.g.~\cite{nocedal2006numerical,bertsekas1982constrained,koehler2022sqp}.

\begin{figure}[tb]
  \small{
  \begin{center}
    \fbox{
      \begin{minipage}{0.95\linewidth}
        \begin{enumerate}[itemsep=0pt, parsep=0pt, topsep=0pt]
        \item[] \textbf{Init:} $\utilde^{(0)}$, $\lambdak[0]$, $\varepsilon_{\text{update}}>0,\varepsilon_{\text{tol}}>0,\mu_{0}> 0$.\vspace{1ex}
        \item[] \textbf{for} $k=0,1,\dots$ until convergence \textbf{do}\vspace*{-0.5ex}
          \begin{enumerate}[label=\arabic*.]
            \item\label{alg:globalized-sqp-1} If $\|{\nabla} \calL^{(k)}\|_{\infty} \leq \varepsilon_{\text{tol}}$, \emph{STOP}.\vspace*{0.5ex}
            \item\label{alg:globalized-sqp-2} Compute the KKT point $(\delta\utilde^{(k)},\delta\lambdak)$ of~\eqref{eq:sqp-subproblem-1}.\vspace*{0.5ex}
            \item\label{alg:globalized-sqp-3} Set $\mu_{k+1} = \max\{\mu_{k},\|\delta\lambdak\|_{\infty}+\varepsilon_{\text{update}}\}$ and compute the step length $\alphak$ which fulfills the Armijo condition~\eqref{eq:penalty-function-armijo} by backtracking.\vspace*{0.5ex}
            \item\label{alg:globalized-sqp-4} Set $\utilde^{k+1} = \utilde + \alphak\delta\utilde^{(k)}$ and $\lambda^{(k+1)} = \delta\lambdak$.
          \end{enumerate}\vspace*{-0.5ex}
        \item[] \textbf{end}
        \end{enumerate}
      \end{minipage}
    }
  \end{center}}
  \caption{Globalized SQP algorithm.}\label{fig:globalized-sqp}
\end{figure}

\subsection{SQP with Quasi-Newton Methods}\label{sec:sqp-with-quasi-newton}
An important feature for the SQP approach is the fact that it does not depend on $\grad^{2}_{\utilde\utilde}\calL$ and hence the Hessian $\grad^{2}_{\utilde\utilde}\calL$ can be replaced by an approximation $H^{(k)}$ which is symmetric and positive definite (s.p.d) on $\ker(B)$, in~\eqref{eq:sqp-subproblem-1} and~\eqref{eq:penalty-function-dir-der}, without losing the theoretical properties,~\cite{nocedal2006numerical,bertsekas1982constrained,koehler2022sqp}.  Quasi-Newton update formulas are an efficient way to calculate an approximation for the Hessian.  In the context of SQP algorithms, for a given approximation $H^{(k)}$, we have to compute the next approximation $H^{(k+1)}$ such that the \emph{secant equation} holds, see, e.g.,~\cite{nocedal2006numerical}.

One of the most efficient Quasi-Newton update formula is the Broyden–Fletcher–Goldfarb–Shanno (BFGS)~\cite{broyden1970convergence,fletcher1970new,goldfarb1970family,shanno1970conditioning} formula.  An important modification for Quasi-Newton methods is the transition to limited-memory variants, where the idea is to use only the last $m$ vectors of $y^{(k)}$ and $d^{(k)}$ for the new approximation.

\begin{figure}[tb]
  \small{
  \begin{center}
    \fbox{
      \begin{minipage}{0.95\linewidth}
        \begin{enumerate}
        \item[] \textbf{Init:} $\utilde^{(0)}$, $\lambdak[0]$, $H^{(0)}$, $\varepsilon_{\text{update}}>0,\varepsilon_{\text{tol}}>0,\mu_{0}> 0$.\vspace{1ex}
        \item[] \textbf{for} $k=0,1,\dots$ until convergence \textbf{do}\vspace*{-0.5ex}
          \begin{enumerate}[label=\arabic*.]
            \item\label{alg:globalized-sqp-with-quasi-newton-1} If $\|{\nabla} \calL^{(k)}\|_{\infty} \leq \varepsilon_{\text{tol}}$, \emph{STOP}.\vspace*{0.5ex}
            \item\label{alg:globalized-sqp-with-quasi-newton-2} Compute the KKT point $(\delta\utilde^{(k)},\delta\lambdak)$ of 
                  \begin{equation}
                    \label{eq:globalized-sqp-with-quasi-newton-approx-problem}
                    \begin{aligned}
                      \textstyle\min_{d}  \grad_{\utilde}{\calL(\utilde^{(k)},\lambdak)}^{T}d + \frac{1}{2}d^{T}H^{(k)}d  \quad\
                      \text{s.t.}\quad  B(\utilde^{(k)} + d) = 0,
                    \end{aligned}
                  \end{equation}
            \item\label{alg:globalized-sqp-with-quasi-newton-3} Set $\mu_{k+1} = \max\{\mu_{k},\|\delta\lambdak\|_{\infty}+\varepsilon_{\text{update}}\}$ and compute the step length $\alphak$ which fulfills the Armijo condition~\eqref{eq:penalty-function-armijo} by backtracking.
            \item\label{alg:globalized-sqp-with-quasi-newton-4} Set $\utilde^{k+1} = \utilde + \alphak\delta\utilde^{(k)}$ and $\lambda^{(k+1)} = \delta\lambdak$.
            \item\label{alg:globalized-sqp-with-quasi-newton-5} Compute $H^{(k+1)}$ according to some Quasi-Newton formula.
          \end{enumerate}\vspace*{-0.5ex}
        \item[] \textbf{end}
        \end{enumerate}
      \end{minipage}
    }
  \end{center}}
  \caption{Globalized SQP algorithm with a full Quasi-Newton approximation.}\label{fig:globalized-sqp-with-quasi-newton}
\end{figure}

An SQP algorithm with a Quasi-Newton approximation is outlined in Fig.~\ref{fig:globalized-sqp-with-quasi-newton}.  The KKT point of~\eqref{eq:globalized-sqp-with-quasi-newton-approx-problem} is know given by the solution of
\begin{equation}
  \label{eq:sqp-kkt-quasi-newton}
  \begin{aligned}
    \begin{pmatrix}
      H^{(k)} &\hspace*{1ex} B^{T} \\[0.5ex]
      B\phantom{^{T}} & 0
    \end{pmatrix}
    \begin{pmatrix}
      \delta\utilde^{(k)} \\[0.5ex]
      \delta\lambdak
    \end{pmatrix}
    = -
    \begin{pmatrix}
      \grad_{\utilde}{\calL(\utilde^{(k)},\lambdak)} \\[0.5ex]
      B\utilde^{(k)}
    \end{pmatrix}
  \end{aligned}.
\end{equation}
Notice, to ensure that we obtain a descent direction for $P_{1}$ by the solution of~\eqref{eq:sqp-kkt-quasi-newton}, it is necessary that $H^{(k)}$ is s.p.d.\ on $\ker(B)$.  The BFGS update formula guarantees that $H^{(k+1)}_{BFGS}$ is s.p.d.\ if $H^{(k)}$ is s.p.d.\ and $\left.y^{(k)}\right.^{T}d^{(k)}$.  In unconstrained optimization, this can be ensured by the Wolfe condition~\cite{wolfe1969convergence}, which is not feasible for $P_{1}$.  Hence, it might happen that $H^{(k+1)}$ is no longer s.p.d.\ and the algorithm fails to compute a descent direction for $P_{1}$.  In this case the whole algorithm breaks.

\subsection{Nonlinear FETI-DP methods with Quasi-Newton}\label{sec:nonlinear-fetidp-with-quasi-newton}

As outlined in~\ref{sec:feti-dp-basics}, nonlinear FETI-DP can be viewed as nonlinear optimization problem with equality constraints.  Thus, we us the SQP algorithm outlined in Fig.~\ref{fig:globalized-sqp-with-quasi-newton} to combine nonlinear FETI-DP with Quasi-Newton methods.  Indeed, we solve~\eqref{eq:minimization-problem} by the SQP algorithm where we approximate $\grad_{\utilde\utilde}\calL$ using the BFGS update formula.

The efficient use the SQP method Fig.~\ref{fig:globalized-sqp-with-quasi-newton}, relies on the fast solution of the subproblem~\eqref{eq:globalized-sqp-with-quasi-newton-approx-problem}, i.e., the fast solution of the KKT system~\eqref{eq:sqp-kkt-quasi-newton}.  Due to the block structure of $\grad_{\utilde\utilde}\calL$ in nonlinear FETI-DP
\begin{equation}
  \label{eq:nl-feti-dp-block}
  \begin{aligned}
    \grad_{\utilde\utilde}\calL(\utilde,\lambda) =
    \begin{pmatrix}
      \grad_{\utilde_{B}\utilde_{B}}\calL(\utilde,\lambda) &\hspace*{1ex}  \grad_{\utilde_{B}\utilde_{\Pi}}\calL(\utilde,\lambda)\\[0.5ex]
      \grad_{\utilde_{\Pi}\utilde_{B}}\calL(\utilde,\lambda) &\hspace*{1ex}  \grad_{\utilde_{\Pi}\utilde_{\Pi}}\calL(\utilde,\lambda)
    \end{pmatrix},
  \end{aligned}
\end{equation}
where the subscript $B$ denote all inner and dual variables and the subscript $\Pi$ denotes the primal variables, a factorization of $\grad_{\utilde\utilde}\calL$ can be computed efficiently.  Notice that the matrix $\grad_{\utilde_{B}\utilde_{B}}\calL(\utilde,\lambda)$ is block diagonal and local to the single subdomains.  Hence, we can use $H^{(0)}:=\grad_{\utilde\utilde}\calL(\utilde^{(0)},\lambda^{(0)})$ and the property of the BFGS approximation that the application of $\left.H^{(k)}\right.^{-1}$ can be efficiently implemented, if $\left.H^{(0)}\right.^{-1}$ is known.  We can solve~\eqref{eq:sqp-kkt-quasi-newton} in the standard FETI-DP way by eliminating $\delta\utilde^{(k)}$ and solve the Schur complement problem
\begin{equation}
  \label{eq:nl-fetidp-quasi-newton-schur-complement}
  \begin{aligned}
    B\left. H^{(k)}\right.^{-1}B^{T}\delta\lambda^{(k)} = g,
  \end{aligned}
\end{equation}
where
\begin{equation*}
  \begin{aligned}
    g = B\utilde^{(k)} - B\left. H^{(k)}\right.^{-1} \grad_{\utilde}{\calL(\utilde^{(k)},\lambdak)},
  \end{aligned}
\end{equation*}
by some Krylov method.  As a preconditioner, we use the Dirichlet preconditioner evaluated at the initial iterate $\utilde^{(0)}$.  Indeed, while we update the BFGS approximation of $\grad_{\utilde\utilde}\calL$, we do \emph{not} update the Dirichlet preconditioner.  In our numerical experiments, this does not increase the number of Krylov iterations for~\eqref{eq:nl-fetidp-quasi-newton-schur-complement} significantly.

Note that it is possible to use some other initial guess for $H^{(0)}$, but in the context of Quasi-Newton methods with application to continuum mechanics, it is crucial to use an exact Hessian as initial value, see, e.g.,~\cite{barnafi2022parallel}.

To accelerate the convergence, we restart our BFGS approximation if some condition of insufficient decrease is fulfilled.  Let $(\utilde^{(k)},\lambda^{(k)})$ be some iterate and $H^{(k)}$ be the corresponding BFGS approximation.  The next iterate $(\utilde^{(k+1)},\lambda^{(k+1)})$ is obtained by the solution of~\eqref{eq:sqp-kkt-quasi-newton} and a step length $\alpha^{(k)}$ which fulfills the Armijo condition~\eqref{eq:penalty-function-armijo}.  If
\begin{subequations}\label{eq:nl-fetidp-quasi-newton-insuf-decrease-1}
  \begin{align}
    | P_{1}(\utilde^{(k+1)};\mu_{k+1}) - P_{1}(\utilde^{(k)};\mu_{k}) | < \eta_{1}|P_{1}(\utilde^{(k)};\mu_{k})|
  \end{align}
  \emph{and}
  \begin{align}
    (1-\eta_{2})\|\grad\calL(\utilde^{(k+1)},\lambda^{(k+1)})\|_{\infty} < \|\grad\calL(\utilde^{(k)},\lambda^{(k)})\|_{\infty},
  \end{align}
\end{subequations}
where $0 < \eta_{1},\eta_{2}<1$, hold, we decide that the search direction, obtained with BFGS approximation $H^{(k)}$, neither reduces the value of the $P_{1}$ nor the first order optimality condition for~\eqref{eq:minimization-problem} enough.  Thus, we compute the exact Hessian and set $H^{(k+1)}:=\grad_{\utilde\utilde}\calL(\utilde^{(k+1)},\lambda^{(k+1)})$ and restart the BFGS approximation.  Note, in this case we have to compute the factorization $\grad_{\utilde\utilde}\calL(\utilde^{(k+1)},\lambda^{(k+1)})$ and also compute the Dirichlet preconditioner for $(\utilde^{(k+1)},\lambda^{(k+1)})$.  If~\eqref{eq:nl-fetidp-quasi-newton-insuf-decrease-1} does \emph{not} hold, we obtain $H^{(k+1)}$ by a BFGS update of $H^{(k)}$.  We outline our nonlinear FETI-DP algorithm with a globalized SQP approach and a Quasi-Newton approximation in~Fig.~\ref{fig:nl-fetidp-globalized-sqp-with-quasi-newton}.

\begin{figure}[tb]
  \small{
  \begin{center}
    \fbox{
      \begin{minipage}{0.95\linewidth}
        \begin{enumerate}
          \item[] \textbf{Init:} $\utilde^{(0)}$, $\lambdak[0]$, $\varepsilon_{\text{update}}>0,\varepsilon_{\text{tol}}>0,\mu_{0}> 0, \eta_{1}\in(0,1),\eta_{2}\in(0,1)$, \par
                \phantom{\textbf{Init:}} $H^{(0)}:=\grad_{\utilde\utilde}\calL(\utilde^{(0)},\lambda^{(0)})$, compute a factorization of $H^{(0)}$.\vspace{1ex}
        \item[] \textbf{for} $k=0,1,\dots$ until convergence \textbf{do}\vspace*{-0.5ex}
          \begin{enumerate}[label=\arabic*.]
            \item\label{alg:nl-fetidp-globalized-sqp-with-quasi-newton-1} If $\|{\nabla} \calL^{(k)}\|_{\infty} \leq \varepsilon_{\text{tol}}$, \emph{STOP}.\vspace*{0.5ex}
            \item\label{alg:nl-fetidp-globalized-sqp-with-quasi-newton-2} Compute the KKT point $(\delta\utilde^{(k)},\delta\lambdak)$ of
                  \begin{equation}
                    \label{eq:nl-fetidp-globalized-sqp-with-quasi-newton-subproblem}
                    \begin{aligned}
                      \min_{d}  \grad_{\utilde}{\calL(\utilde^{(k)},\lambdak)}^{T}d + \frac{1}{2}d^{T}H^{(k)}d  \quad\
                      \text{s.t.}\quad  B(\utilde^{(k)} + d) = 0,
                    \end{aligned}
                  \end{equation}
            \item\label{alg:nl-fetidp-globalized-sqp-with-quasi-newton-3} Set $\mu_{k+1} = \max\{\mu_{k},\|\delta\lambdak\|_{\infty}+\varepsilon_{\text{update}}\}$ and compute the step length $\alphak$ which fulfills the Armijo condition~\eqref{eq:penalty-function-armijo} by backtracking.\vspace*{0.5ex}
            \item\label{alg:nl-fetidp-globalized-sqp-with-quasi-newton-4} Set $\utilde^{k+1} = \utilde + \alphak\delta\utilde^{(k)}$ and $\lambda^{(k+1)} = \delta\lambdak$.\vspace*{0.5ex}
            \item\label{alg:nl-fetidp-globalized-sqp-with-quasi-newton-5} \textbf{If}~\eqref{eq:nl-fetidp-quasi-newton-insuf-decrease-1} holds, \textbf{then}\par
                  \hspace*{4ex} Set $H^{(k+1)}:=\grad_{\utilde\utilde}\calL(\utilde^{(k+1)},\lambda^{(k+1)})$ \par
                  \hspace*{4ex} and compute a factorization of $H^{(k+1)}$\par
                  \textbf{else}\par
                  \hspace*{4ex} Compute $H^{(k+1)}$ by the BFGS update.
          \end{enumerate}\vspace*{-0.5ex}
        \item[] \textbf{end}
        \end{enumerate}
      \end{minipage}
    }
  \end{center}}
  \caption{Nonlinear FETI-DP with a globalized SQP algorithm and a Quasi-Newton approximation.}\label{fig:nl-fetidp-globalized-sqp-with-quasi-newton}
\end{figure}

Let us remark that there is a trade-off between the classical SQP algorithm~Fig.~\ref{fig:globalized-sqp} for nonlinear FETI-DP and the SQP algorithm~Fig.~\ref{fig:nl-fetidp-globalized-sqp-with-quasi-newton} with a Quasi-Newton approximation, due to the factorization of $\grad_{\utilde\utilde}\calL$ within nonlinear FETI-DP.  Normally, the classical SQP algorithm needs a fewer number of SQP subproblems to converge, due to the usage of the exact Hessian, but although the factorization of $\grad_{\utilde\utilde}\calL$ can be implemented very efficiently in parallel, there is a potential bottleneck in the (local) factorization of $\grad_{\utilde_{B}\utilde_{B}}\calL$ and also the factorization of the Schur complement of $S_{\Pi\Pi}=\grad_{\utilde_{\Pi}\utilde_{\Pi}}\calL - \grad_{\utilde_{\Pi}\utilde_{B}}\calL\grad_{\utilde_{B}\utilde_{B}}\calL^{-1}\grad_{\utilde_{B}\utilde_{\Pi}}\calL$ due to the runtime complexity of the sparse direct solvers.  On the other hand, we can reduce the number of factorization of $\grad_{\utilde\utilde}\calL$ with our Quasi-Newton approach, but we need more SQP subproblems to converge and the solution of a subproblem relies on some Krylov method.  Hence, we increase the total amount of Krylov iterations, which also increases the communication in a parallel implementation.

\section{Numerical Results}\label{sec:numerical-results}

We consider a two- or three-dimensional beam bending benchmark problem with a compressible Neo-Hookean constitutive law using $Q$$2$ finite elements.  The strain energy density function is given by $J(x) = \frac{\mu}{2}(\text{tr}(F(x)^{T}F(x)) - 3) - \mu\log(\psi(x)) + \frac{\lambda}{2}(\log(\psi(x)))^{2}$, where $\psi(x) = \det(F(x))$, $F(x) = {\nabla} \varphi(x)$, $\varphi(x)= x + u(x)$, $u(x)$ denotes the displacement and $\mu$ and $\lambda$ are the Lamé constants.  As material parameters, we use $E = 210$ and $\nu = 0.3$.

As a baseline, we use a globalized Newton-like algorithm based on the exact differentiable penalty function $P(\utilde,\lambda;\mu) := \calL(\utilde,\lambda) + \frac{\mu}{2}\|B\utilde\|^{2} + \|B\grad_{\utilde}\calL(\utilde,\lambda)\|$ introduced by Di~Pillo and Grippo in \cite{di1979new}.  Note that a Newton-like direction for $P$ can be computed by the solution of the standard Lagrange-Newton equation.  For the combination of nonlinear FETI-DP methods and the penaly function $P$ see, e.g.,\cite{koehler2022globalization,koehler2022globalization:phd}.

\begin{table}[h]
  \scriptsize{
  \caption{Comparison of Quasi-Newton SQP approach with a Newton-like line search method for $P$ in 2D, local subdomains = $160\times 160$, \#dofs subdomain = 206\,082}\label{table:sqp-penalty-2d}
  \centering
  \begin{tabularx}{0.6\textwidth}{ r r r r }
    \toprule
    \multicolumn{4}{c}{Sizes for local subdomains of $160\times 160$ $Q_{2}$-elements}  \\
    \midrule
    \#Subds  & \#Dofs ($\tilde{x},\lambda$) & \#Dofs ($\tilde{x}$) & \#Coarse Dofs \\
    \midrule
      40     &    8\,280\,168               &     8\,242\,856      &       348     \\[0.5ex]
     160     &   33\,148\,812               &    32\,971\,152      &    1\,500     \\[0.5ex]
     360     &   74\,605\,936               &    74\,184\,888      &    3\,452     \\[0.5ex]
     640     &  132\,651\,540               &   131\,884\,064      &    6\,204     \\[0.5ex]
     \bottomrule
  \end{tabularx}
  \vspace*{3ex}

   \begin{tabularx}{1.\textwidth}{ r | r r r r | r r r}
    \toprule
     &   \multicolumn{4}{c}{Quasi-Newton SQP}                               &   \multicolumn{3}{c}{Newton-like method for $P$}\\
     \#Subds  & Solve & \#It. & \#Krylov\,It.  & \#Recomp ($\grad^{2}_{\utilde\utilde}\calL$)   &  Solve & \#Newton\,It. & \#Krylov\,It.      \\
    \midrule
      40     &      321.1s   &   25   &      589       &             8           &     400.7s   &  15             &  355  \\[0.5ex]
     160     &      421.0s   &   26   &      624       &             8           &     514.8s   &  14             &  341  \\[0.5ex]
     360     &      381.0s   &   24   &      578       &             7           &     522.8s   &  14             &  349  \\[0.5ex]
     640     &      422.4s   &   25   &      615       &             8           &     525.7s   &  14             &  355  \\[0.5ex]
     \bottomrule
   \end{tabularx}
 }
\end{table}

\begin{table}[h]
  \scriptsize{
    \caption{Comparison of Quasi-Newton SQP approach with a Newton-like line search method for $P$ in 3D, local subdomains = $18\times 18\times 18$, \#dofs subdomain = 151\,959}\label{table:sqp-penalty-3d}
  \centering
  \begin{tabularx}{0.6\textwidth}{ r r r r }
    \toprule
    \multicolumn{4}{c}{Sizes for local subdomains of $18\times 18\times 18$ $Q_{2}$-elements}  \\
    \midrule
    \#Subds  & \#Dofs ($\tilde{x},\lambda$) & \#Dofs ($\tilde{x}$) & \#Coarse Dofs \\
    \midrule
      80     &   12\,815\,196               &    12\,152\,964      &    1\,800             \\[0.5ex]
     640     &  104\,010\,420               &    97\,213\,008      &   16\,020             \\[0.5ex]
     \bottomrule
   \end{tabularx}

   \vspace*{3ex}

   \begin{tabularx}{1.\textwidth}{ r | r r r r | r r r}
     \toprule
     &   \multicolumn{4}{c}{Quasi-Newton SQP}                               &   \multicolumn{3}{c}{Newton-like method for $P$}\\
     \#Subds  & Solve & \#It. & \#Krylov\,It.  & \#Recomp\,($\grad^{2}_{\utilde\utilde}\calL$)   &  Solve & \#Newton\,It. & \#Krylov\,It.      \\
     \midrule
      80     &   939.0s     &   22   &    1\,415     &   5                 &   2\,010.9s  &  14             &       948    \\[0.5ex]
     640     &   486.9s     &   11   &       596     &   2                 &   1\,818.1s  &  10             &       857    \\[0.5ex]
     \bottomrule
   \end{tabularx}

   }
\end{table}

In Table~\ref{table:sqp-penalty-2d}, we show weak scaling results for the two-dimensional problem comparing our SQP approach using Quasi-Newton approximation of the Hessian with the base line using the exact differentiable penalty function $P$.  We increase the number of subdomains from $40$ up to $640$ using subdomains with $206\,082$ degrees of freedom.  We report the nonlinear solution time, the number of (nonlinear) iterations, and the total number of Krylov iterations.  Note that for both methods one nonlinear iteration involves the solution of the saddle point problem \eqref{eq:sqp-kkt-quasi-newton}, where we use a Quasi-Newton approximation for the SQP approach and the exact Hessian for the penalty function $P$.  But for the Newton-like method based on $P$, we need to compute a factorization of $\grad^{2}_{\utilde\utilde}\calL$ in each nonlinear iteration which is not the case for the SQP approach using Quasi-Newton approximations, where we only need to compute such a factorization if we recompute the original Hessian of $\calL$ as outlined in Fig.~\ref{fig:nl-fetidp-globalized-sqp-with-quasi-newton}.  Hence, the computational saving is in the fewer factorizations of $\grad^{2}_{\utilde\utilde}\calL$.  We need between 7 and 8 factorizations for the SQP approach, whereas we need between 14 and 15 factorizations for $P$, which means we can save nearly 50\% of the factorizations.  Nevertheless, we need more (nonlinear) iterations for the SQP approach, between 24 and 26, compared to the 14 or 15 iterations for $P$, resulting in a higher number of (total) Krylov iterations for the SQP approach.  The higher number of Krylov iterations consumes some of the time savings.  Overall, we obtain a speed up of the solution between 18\% for 160 subdomains and 27\% for 360 subdomains for the Quasi-Newton method using SQP approach against our base line based on $P$.

In Table~\ref{table:sqp-penalty-3d}, we show weak scaling results for the three-dimensional system.  Due to the structured decomposition in our experiments, we only report results for 80 and 640 subdomains.  The number of factorizations of $\grad^{2}_{\utilde\utilde}\calL$ in the SQP approach reduces drastically from 14 and 10 for the penalty method based on $P$ to 5 and 2, for 80 and 640 subdomains, respectively.  Therefore, also the nonlinear solution times reduces over 50\%.  This shows that Quasi-Newton methods can accelerate the solution time by reducing the number of factorizations of $\grad^{2}_{\utilde\utilde}\calL$.

Note that using a Quasi-Newton method for $P$ is very challenging, since $\grad^{2}\calL$ is not a Newton-like approximation for $\grad^{2}P$ only the solution of the standard Lagrange-Newton equation reveals a Newton-like search direction.  Furthermore, the recomputation of $\grad^{2}_{\utilde\utilde}\calL$ is essential to obtain a good speed up, i.e., without the recomputation the SQP approach with Quasi-Newton approximations is slower than the Newton-like method for $P$.

\bibliographystyle{spmpsci}
\bibliography{./all_refs}

\end{document}